\documentclass[12pt]{article}
\newcommand{\fr}{ \frac}
\newcommand{\beq}{\begin{equation}}
\newcommand{\eeq}{\end{equation}}
\newcommand{\br}{\begin{array}}
\newcommand{\er}{\end{array}}
\newcommand{\lb}{\label}
\newcommand{\ar}{\rightarrow}

\begin{document}
 \begin{center}
 {\Large\bf  Fractional Super Lie Algebras and Groups }
 \end{center}
 \vspace{5mm}
 \begin{center}
 H. Ahmedov$^1$, A. Yildiz$^1$ and Y. Ucan$^2$
 \end{center}

 \noindent
 1. Feza G\"ursey Institute, P.O. Box 6, 81220,
 \c{C}engelk\"{o}y, Istanbul, Turkey
 \footnote{E-mail : hagi@gursey.gov.tr }. \\
 2. Yildiz Technical University, Department of Mathematics,
 Besiktas, Istanbul, Turkey.
 \vspace{5mm}
 \begin{center}
 {\bf Abstract}
 \end{center}
 $\mathrm{n^{\underline{th}}}$ root of a Lie algebra and its dual (that is fractional
 supergroup )  based on  the permutation group $S_n$ invariant forms is formulated
 in the Hopf algebra formalism. Detailed discussion of   $S_3$-graded $\mathrm{sl(2)}$
 algebras is done.

 \noindent

 \vspace{1cm}
 \noindent
 {\bf 1. Introduction}

 \vspace{5mm}
 \noindent
 To arrive at a superalgebra  one adds new elements $Q_\alpha$ to  generators
 $X_j$ of the corresponding Lie algebra and define relations
 \beq\lb{z2}
 \{Q_\alpha, Q_\beta\}=b^j_{\alpha\beta}X_j.
 \eeq
 Observing that the anticommutator in the above relation is invariant
 under the cyclic $\mathrm{Z}_2$ or permutation $\mathrm{S}_2$ groups  we can look
 for possible generalization of the supersymmetry by using
 $S_n$ or $Z_n$ invariant structures instead of the anticommutator.
 For example if   $n=3$ instead of  (\ref{z2}) one has cubic relations
 \beq\lb{z3}
 Q_\alpha Q_\beta Q_\gamma + Q_\gamma Q_\alpha Q_\beta
 + Q_\beta Q_\gamma Q_\alpha =
 b^j_{\alpha\beta\gamma}X_j
 \eeq
 which is $\mathrm{Z_3}$ invariant and
 \beq\lb{s3}
 Q_\alpha\{Q_\beta, Q_\gamma\} +Q_\beta\{Q_\alpha,
 Q_\gamma\}+
 Q_\gamma\{Q_\alpha, Q_\beta\}=b^j_{\alpha\beta\gamma}X_j
 \eeq
 which is $\mathrm{S_3}$ invariant.
 From the above relations only (\ref{s3}) appears to be consistent at the
 co-algebra level.  Sometimes we will use the term fractional
 superalgebras for  $\mathrm{S_n}$-graded algebras with $n=3,4,\dots ...$
 with fractional super groups  being their dual.

 Fractional super algebras based on  $S_n$ invariant forms  were
 first introduced in \cite{rau1,rau2}. In the present paper we put this construction
 in the Hopf algebra context and define their dual, that is fractional supergroups.
 There are many reasons for doing that. In the formulation of superalgebras one can
 use either geometric   (See for example \cite{dewitt}) or algebraic \cite{kos} approaches
 (See  also  \cite{bat} for a comparison ). As for  fractional superalgebras geometric approach
 seems to be insufficient. This situation is similar to the theory of quantum
 algebras, where we have to work with universal enveloping algebras rather
 than with Lie algebras \cite{vil1}. Moreover, having a fractional
 superalgebra in hand we can define fractional supergroups by taking
 the dual of the former. And at last having put fractional
 superalgebras in the Hopf algebra context we can use well
 developed representation theory of the latter in the construction of
 representations of fractional superalgebras.

 There are another approaches to fractional supersymmetry in literature
 \cite{ahn,rau,ker,vik,azc,dur}.
 For example one can arrive at fractional supergroups by using
 quantum groups at roots of unity \cite{ahm1}.

 The plan of the paper is as follows.  To make the treatment reasonably self-consistent,
 in  section 2 we  give a formulation of  super algebras and
 groups in the Hopf algebra formalism.
 In section 3 we define fractional super algebras and discuss the
 structure of their dual ( fractional super groups ).
 Section 4 is devoted to the detailed discussion of $S_3$-graded  $\mathrm{sl(2)}$
 algebras.

 \vspace{5mm}
 \noindent

 \vspace{1cm}
 \noindent
 {\bf 2. Preliminaries on super algebras }

 \vspace{2mm}
 \noindent
  Let $\mathrm{U(g)}$ be the universal enveloping algebra of a Lie algebra $\mathrm{g}$
 generated by $X_j$, $j=1,\dots \mathrm{dim (g)}$ with
 \beq\lb{la}
 [X_i, X_j]=\sum_{k=1}^{\mathrm{dim(g)}}c^k_{ij}X_k,
 \eeq
 where $c^k_{ij}$ are the structure constants of the Lie algebra  $\mathrm{g}$.
 The Hopf algebra structure of $\mathrm{U(g)}$ is given by the
 co-multiplication  $\Delta:\mathrm{U(g)\ar U(g)\otimes U(g)}$,
 co-unit  $\varepsilon:\mathrm{ U(g)\ar C}$  and antipode  $S: \mathrm{U(g)\ar
 U(g)}$:
 \beq\lb{cla}
 \Delta(X_j)=X_j\otimes 1 +1\otimes X_j, \ \ \ \varepsilon(X_j)=0, \ \ \ S(X_j)=-X_j
 \eeq
 We can extend the Hopf algebra $\mathrm{U(g)}$ by adding
 elements $Q_\alpha$, $\alpha =1, \dots, N$ and $K$ with relations
 \beq\lb{s2al1}
 \{Q_\alpha,Q_\beta\}= \sum_{j=1}^{\mathrm{dim(g)}}b_{\alpha\beta}^jX_j
 \eeq
 \beq\lb{s2al2}
 [Q_\alpha, X_j]=
 \sum_{\beta=1}^{\mathrm{N}}a^j_{\alpha\beta}Q_\beta,
 \eeq
 \beq
 KQ_\alpha =-Q_\alpha K, \ \ \ \ K^2=1
 \eeq
 where $b_{\alpha\beta}^j$ and $a^j_{\alpha\beta}$ are the structure
 coefficients  satisfying the super Jacobi identities.
 This algebra which  we denote by $\mathrm{U_2^N(g)}$ can also be equipped with a Hopf
 algebra structure  by defining
 \beq
 \Delta (Q_\alpha)= Q_\alpha\otimes 1+ K\otimes Q_\alpha, \ \ \Delta (K) = K\otimes K, \ \
 \eeq
 \beq
 \varepsilon(Q_j)=0, \ \ \varepsilon(K)=1, \ \ S(Q_j)=Q_jK, \ \  S(K)=K.
 \eeq

 The dual of $\mathrm{U_2^N(g)}$ is the Hopf algebra $\mathrm{A_2^N(G)=
 C^\infty(G)\times\Lambda_2^N}$, where  $\mathrm{C^\infty(G)}$ is the algebra of infinite
 differentiable  functions on a Lie group $\mathrm{G}$ and  $\Lambda_2^N$
 is the algebra over the field of complex numbers  generated  by
 $\theta_\alpha$, $j=1,\dots,N$ and  $\lambda$ with relations
 \beq
 \{\theta_\alpha,\theta_\beta\} = 0, \ \ \ \{\lambda, \theta_\alpha\}=0, \ \ \  \lambda^2 =1.
 \eeq
 The operations $\Delta$, $\varepsilon$ and $S$ in $\mathrm{A^N(G)}$ depend on the value of
 structure constants $c^k_{ij}$, $b_{\alpha\beta}^j$ and $a^j_{\alpha\beta}$.

 For example if $\mathrm{G=C^N}$ then the formulas
 \beq
 \Delta(\theta_\alpha)=\theta_\alpha\otimes 1 + \lambda\otimes \theta_\alpha, \ \
 \Delta(\lambda)=\lambda\otimes\lambda,\ \
 \Delta(z_\alpha)= z_\alpha\otimes 1 + 1\otimes
 z_\alpha+\lambda\theta_\alpha\otimes\theta_\alpha,
 \eeq
 \beq
 \varepsilon(\theta_\alpha)=0, \ \ \ \varepsilon(\lambda)=1, \ \ \ \varepsilon(z_\alpha)=0
 \eeq
 \beq
  S(\theta_\alpha)=-\lambda\theta_\alpha, \ \ \ S(\lambda)=\lambda, \ \ \
  S(z_\alpha)=-z_\alpha
 \eeq
 define the super $\mathrm{N}$-dimensional translation group.
 The corresponding  super algebra is defined by
 \beq
 \{Q_\beta,Q_\alpha\}=\delta_{\beta\alpha}P_\alpha, \ \ \  [X_\beta, X_\alpha]=0,
 \ \ \ [Q_\beta, X_\alpha]=0.
 \eeq

\vspace{5mm}
 \noindent
 {\bf 3. Fractional superalgebras and supergroups}

 \vspace{2mm}
 \noindent
 To arrive at cubic root of a lie algebra $g$ we have to replace $S_2$
 invariant form in (\ref{s2al1}) by $S_3$ invariant one. Consequently we
 define an algebra generated by $X_j$, $j=1,\dots, dim(g)$ and
 $Q_\alpha$, $K$, $\alpha =1,\dots, N$ satisfying the relations
 (\ref{la}) and
 \beq\lb{s3al1}
 \{ Q_\alpha,Q_\beta,Q_\gamma\}= b^j_{\alpha\beta\gamma}X_j,
 \eeq
 \beq\lb{s3al2}
 [Q_\alpha, X_j]= a^j_{\alpha\beta} Q_\beta,
 \eeq
 and
 \beq\lb{s3al3}
 KQ_\alpha =qQ_\alpha K, \ \ \ q^3=1, \ \ \ K^3=1,
 \eeq
 where
 \beq
 \{ Q_\alpha, Q_\beta, Q_\gamma \} \equiv
 Q_\alpha \{Q_\beta, Q_\gamma\} + Q_\beta \{Q_\alpha, Q_\gamma\}+
 Q_\gamma \{Q_\alpha, Q_\beta\}
 \eeq
 is $S_3$ invariant form. We denote this algebra by the symbol
 $U^N_3(g)$ with the lower index indicating the degree of
 grading. One can check that the above algebra is compatible with
 the co-algebra structure and antipode given by the formulas
 \beq\lb{com}
 \Delta (Q_\alpha)= Q_\alpha\otimes 1+ K\otimes Q_\alpha, \ \ \Delta (K) = K\otimes K, \ \
 \eeq
 \beq\lb{ant}
 \varepsilon(Q_j)=0, \ \ \varepsilon(K)=1, \ \ S(Q_j)= -K^2Q_j, \ \  S(K)=K^2.
 \eeq
 For example let us verify the consistency of the comultiplication
 $\Delta$ with (\ref{s3al1}). Since $\Delta$ is a homomorphism we
 have
 \begin{eqnarray}
 \Delta (Q_\alpha Q_\beta Q_\gamma)=Q_\alpha Q_\beta Q_\gamma\otimes 1
 + 1\otimes Q_\alpha Q_\beta Q_\gamma +
 Q_\alpha Q_\beta K\otimes  Q_\gamma  \nonumber \\
 +  Q_\alpha K Q_\gamma\otimes  Q_\beta +  KQ_\beta Q_\gamma\otimes  Q_\alpha +
 Q_\alpha K^2\otimes Q_\beta Q_\gamma \nonumber \\
 + KQ_\beta K\otimes Q_\alpha Q_\gamma +
 K^2Q_\gamma \otimes Q_\alpha Q_\beta
 \end{eqnarray}
 Using (\ref{s3al3}) we get
 \beq
 \sum_{(\alpha\beta\gamma)\in S_3}(
 Q_\alpha Q_\beta K\otimes  Q_\gamma +
 Q_\alpha K Q_\gamma\otimes  Q_\beta
 + KQ_\beta Q_\gamma\otimes  Q_\alpha )= 0
 \eeq
 and
 \beq
 \sum_{(\alpha\beta\gamma)\in S_3}(
 Q_\alpha K^2\otimes Q_\beta Q_\gamma + KQ_\beta K\otimes Q_\alpha Q_\gamma +
 K^2Q_\gamma \otimes Q_\alpha Q_\beta )=0.
 \eeq
 Thus we have shown that
 \beq
 \sum_{(\alpha\beta\gamma)\in S_3} \Delta (Q_\alpha Q_\beta
 Q_\gamma)=
 \sum_{(\alpha\beta\gamma)\in S_3}( Q_\alpha Q_\beta
 Q_\gamma\otimes 1+ 1\otimes  Q_\alpha Q_\beta
 Q_\gamma)
  \eeq
 which together with the comultiplication rule (\ref{cla}) for the generators
 $X_j$ implies the consistency of the comultiplication  (\ref{com}) with the
 relation (\ref{s3al1}).

 To define structure constants $b^j_{\alpha\beta\gamma}$ and
 $a^j_{\alpha\beta}$ we have to derive identities involving
 commutator and $S_3$ invariant form. One can check that relations
 \beq\lb{j1}
 [A, [B, C]]+[C, [A, B]]+[B, [C, A]]=0,
 \eeq
 \beq\lb{j2}
 [A, \{B, C, D\}]+\{[B, A], C, D\}+\{B, [C, A], D\}+\{B, C, [D,
 A]\} =0
 \eeq
 and
 \beq\lb{j3}
 [A, \{B, C, D\}]+ [B, \{A, C, D\}]+ [C, \{B, A, D\}]+[D, \{B, C, A\}]=0
 \eeq
 are satisfied identically. For example we verify the identity
 (\ref{j2}). Let $B=A_1$, $C=A_2$ and $D=A_3$. Then
 \beq
 [A,\{A_1, A_2, A_3\}]=\sum_{ijk\in S_3} ( [A,A_i] A_jA_k +
 A_iA_j[A,A_k] + A_i[A,A_j]A_k ).
 \eeq
 Combining terms $(123), \ (132)$ from the first sum on the right hand side of
 the above equality, $(231), \ (321)$ from the second sum and $(213), \ (312)$
 from the third sum we get $\{[A,A_1], A_2, A_3\}$. In a similar
 fashion we obtain $\{[A,A_2], A_1, A_3\}$ and $\{[A,A_3], A_2,
 A_1\}$. Thus
 \beq
 [A,\{A_1, A_2, A_3\}]=\{[A,A_1], A_2, A_3\}+ \{[A,A_2], A_1,
 A_3\}+  \{[A,A_3], A_2, A_1\}
 \eeq
 which is the identity (\ref{j2}).

 The one given by (\ref{j1}) is the  usual Jacobi identity. Inserting
 \beq
 A=X_i, \ \ \ B=X_j, \ \ \ C=Q_\alpha
 \eeq
 into (\ref{j1}) and  using  (\ref{s3al2}) and (\ref{la}) we get
 \beq\lb{rep1}
 \sum_{\sigma=1}^N (a^i_{\alpha\sigma}a^j_{\sigma\beta}-
 a^j_{\alpha\sigma}a^i_{\sigma\beta})= \sum_{k=1}^{dim (g)}
 c^k_{ij}a^k_{\alpha\beta}.
 \eeq
 Comparing the above relation with (\ref{la}) we conclude that the
 $N\times N$  matrices $a^j\equiv (a^j_{\alpha\beta})_{\alpha,\beta=1}^N$ define
 an $N$-dimensional  representation of a given Lie algebra.
 There are different possibilities in the choice of this
 representation. For example if $g=sl(2)$ and $N=2$ we can either use the
 scalar representation  $a^j_{\alpha\beta}=0$ or the spinor one in which
 $a^j_{\alpha\beta}$ are the Pauli matrices.  Consequently for fixed Lie algebra $g$ and $N$ we can define  different
 super fractional  algebras. To be more precise one has to add
 additional index in the notation $U_3^N(g)$ which reflects
 the transformation law of super generators $Q_\alpha$ with respect to a given Lie algebra
 $g$. However for the sake of simplicity we will not do it.
 Detailed discussion of this non uniqueness is done  in the next section where we
 consider fractional super algebras $sl(2)$.

 Let us now consider restrictions on structure coefficients coming from the other identities.
 Inserting
 \beq
 A=X_k, \ \ \ B=Q_\alpha, \ \ \ C=Q_\beta, \ \ \ D=Q_\gamma
 \eeq
 into the identity (\ref{j2}) and
 \beq
 A=Q_\sigma, \ \ \ B=Q_\alpha, \ \ \ C=Q_\beta, \ \ \ D=Q_\gamma
 \eeq
 into (\ref{j3}) and using (\ref{s3al1}),  (\ref{s3al2}) we arrive
 at relations
 \beq\lb{s3j1}
 \sum_{\sigma=1}^N (a^k_{\alpha\sigma}b_{\sigma\beta\gamma}^i +
 a^k_{\beta\sigma}b_{\sigma\alpha\gamma}^i + a^k_{\gamma\sigma} b^i_{\sigma\beta\alpha})=
 \sum_{j=1}^{dim g}c_{jk}^ib^{j}_{\alpha\beta\gamma}
 \eeq
 and
 \beq\lb{s3j2}
 \sum_{k=1}^{dim g} (b^k_{\alpha\beta\gamma}a^k_{\sigma\tau}+
 b^k_{\sigma\alpha\beta}a^k_{\gamma\tau}+b^k_{\gamma\sigma\alpha}
 a^k_{\beta\tau}+ b^k_{\beta\gamma\sigma}a^k_{\alpha\tau})=0.
 \eeq
 Now we define fractional super groups.
 Let $x=\{x_{nm}\}$ be the matrix  representing  a Lie group $G$ and $A(G)$ be the algebra
 of polynomials on $G$. It is known that $A(G)$ is the Hopf
 algebra which is in non degenerate duality with the universal
 enveloping algebra $U(g)$ \cite{vil}. In general the number of group elements
 $x_{nm}$ is more than the number of generators $X_j$ in the corresponding
 Lie algebra $g$.  This is due to the fact that there may be some restrictions on
 the matrix representing a Lie group.
 For example if $G=SL(2)$ we have two by two matrix with
 determinant equal to 1. The number of  independent group parameters
 is equal to the number of generators of $sl(2)$. Explicitly one
 can define these parameters by using some decomposition ( Gauss,
 Cartan, Iwasawa and so on ). In a similar way for an arbitrary
 matrix Lie group we can resolve restrictions imposed on the
 elements $x_{nm}$ and obtain  independent  group parameters $x_j$
 with duality relations
 \beq
 \langle x_i, X_j\rangle =\delta_{ij},
 \eeq
 where $X_j$ are  the generators  of the corresponding Lie algebra.
 However in general $A(G)$ in terms of these new parameters will
 not be the polynomial algebra. It appears that in the Hopf
 algebra formalism it is more convenient to work with elements
 $x_{nm}$. Instead of solving restrictions imposed on these
 elements
 one  defines new generators $X_{nm}$ with some restrictions.
 For example if $g=sl(2)$ we define four generators
 with the restriction $X_{11}+X_{22}=0$.

 To  construct the  dual  algebra to a fractional super algebra
 $U_3^N(g)$ we have to introduce new parameters $\theta_\alpha$,
 $\alpha=1,\dots, N$ and $\lambda$ corresponding to the
 fractional super generators $Q_\alpha$ and $K$. The duality
 relations are given by the following formulas
 \beq\lb{pair}
 \langle\theta_{\alpha}, Q_{\beta}\rangle =\delta_{\alpha\beta}, \
 \ \ \ \langle\lambda, K\rangle =q, \ \ \ \langle x_{nm}, K\rangle =\delta_{nm}
 \eeq
 with all other linear relations being zero.
 Recall the property of the duality relations \cite{vil}
 \beq\lb{dual}
 \langle ab,\phi\rangle=\sum_j\langle a,\phi_j\rangle\langle b,\phi^\prime_j\rangle
 \eeq
 with
 \beq
 \Delta(\phi)=\sum_j \phi_j\otimes \phi^\prime_j.
 \eeq
 Here $\phi$ and $a, \ b$ are elements of a Hopf algebra and its
 dual. Inserting in (\ref{dual}) $a=\theta_\alpha$, $b=\lambda$ and $\phi=Q_\alpha$
 and using (\ref{com}), (\ref{pair}) we get
 \beq\lb{f1}
 \lambda\theta_\alpha =q\theta_\alpha \lambda
 \eeq
 Taking $a=x_{nm}$, $b=\lambda$ and $\phi=X_{nm}$ we
 conclude that elements $x_{nm}$ commute with $\lambda$. The
 choice $a=\lambda^2$, $b=\lambda$ and $\phi = K$ implies
 $ \langle \lambda^3, K\rangle =1$. Since $\lambda^3$ cannot be proportional to the diagonal
 elements $x_{nn}$ ( $\langle\lambda^3, X_{nn}\rangle=0$ ) we  have
 \beq\lb{f2}
 \lambda^3=1.
 \eeq
 The above condition can be shown to imply the  comultiplication
 \beq\lb{f3}
 \Delta (\lambda )= \lambda\otimes \lambda.
 \eeq
 To make (\ref{f1}) and  (\ref{f3}) compatible we have to define
 \beq\lb{f4}
 \Delta (\theta_\alpha)=\sum_{\beta=1}^N\theta_\beta\otimes d_{\beta\alpha}+
 \lambda \otimes \theta_{\beta} + \cdots,
 \eeq
 where $d=\{d_{\alpha\beta}\}$ is a N dimensional representation of a Lie group $G$
 under consideration and  $\cdots$ denotes the
 combination of terms consisting of $4$, $7$, $10$ and so on generators
 $\theta_\alpha$. Using (\ref{f4}), (\ref{pair}) and (\ref{dual}) after
 some algebra we get
 \beq\lb{f6}
 \{\theta _\alpha, \theta_\beta, \theta_\gamma\}=0
 \eeq
 In a similar way  the commutativity of $x_{nm}$ with
 $\lambda$ and (\ref{f3}) imply
 \beq\lb{f5}
 \Delta (x_{nm})= \sum_k x_{nk}\otimes x_{km} + \cdots,
 \eeq
 where $\cdots$ denotes the
 combination of terms consisting of $3$, $6$, $9$ and so on generators
 $\theta_\alpha$. Using (\ref{f5}) and  (\ref{dual}) we conclude that
 elements $x_{nm}$ commute with $\theta_\alpha$.

 Let us denote the algebra generated by $\theta_\alpha$,
 $\alpha=1, \dots N$ and $\lambda$ satisfying (\ref{f1}), (\ref{f2}) and (\ref{f6})
 by $\Lambda_3^N$ and the
 direct product algebra $A(G)\times \Lambda^N_3$ by $A_3^N(G)$.
 This algebra   is in  nondegenerate duality with a
 Hopf algebra $U_3^N(g)$.  We call $A_3^N(G)$ fractional supergroup.
 Using  the properties
 \beq
 \varepsilon (a)=\langle a, 1\rangle
 \eeq
 and
 \beq
 \langle S(a), \phi \rangle= \langle a, S(\phi )\rangle
 \eeq
 of the duality relations we get the counit operation
 \beq
 \varepsilon (x_{nm})=\delta_{nm}, \ \  \
 \varepsilon (\theta_\alpha )=0, \ \ \  \varepsilon (\lambda)=1.
 \eeq
 and the antipode
 \beq
 S(\lambda )= \lambda^2.
 \eeq
 Using the properties of duality relations and  Hopf algebra axioms one can
 derive  unknown terms in (\ref{f4}) and (\ref{f5}) and antipodes $S(x_{nm})$,
 $S(\theta_{\alpha})$. These calculations depend on structure constants $c^i_{jk}$,
 $a^j_{\alpha\beta}$ and $b^j_{\alpha\beta\gamma}$. We demonstrate this construction on the
 explicit examples which will be given later.

 Before closing this section we define  $S_n$-graded Lie algebras and groups.
 This can be done in the same way as $S_3$ case.
 For this one has to use  $S_n$ invariant form
 \beq
 \{Q_{\alpha_1}, Q_{\alpha_2},\dots, Q_{\alpha_n}\}=
 \sum_{\alpha_1,\alpha_2,\dots,\alpha_n \in S_n}Q_{\alpha_1}Q_{\alpha_2}\cdots
 Q_{\alpha_n},
 \eeq
 where summation runs over all permutations of $S_n$. Instead of
 (\ref{s3al1}) and (\ref{s3al3}) we then have
 \beq
 \{Q_{\alpha_1}, Q_{\alpha_2},\dots,
 Q_{\alpha_n}\}=b^j_{\alpha_1\alpha_2,\dots,\alpha_n} X_j
 \eeq
 and
 \beq
 KQ_\alpha =qQ_\alpha K, \ \ \ q^n=1, \ \ \ K^n=1
 \eeq
 such that
 \beq\lb{snj1}
 \sum_{\sigma=1}^N\sum_{(\alpha_1,\dots, \alpha_n)\in Z_n}
 a^k_{\alpha_1\sigma}b_{\sigma_2\dots\alpha_n}^i =
 \sum_{j=1}^{dim g}c_{jk}^ib^{j}_{\alpha_1\alpha_2\dots\alpha_n}
 \eeq
 and
 \beq\lb{snj2}
 \sum_{k=1}^{dim(g)}\sum_{(\alpha_1,\dots,\alpha_{n+1})\in Z_{n+1}}
 b^k_{\alpha_1\dots \alpha_n}a^k_{\alpha_{n+1}\tau}=0.
 \eeq
 The multiplication and counit  in $U^N_n(g)$ are similar to that
 in $U^N_2(g)$ or $U^N_3(g)$ while the antipode is given by
 \beq
 S(Q_\alpha)=-K^{n-1}Q_\alpha, \ \ \ S(K)=K^{n-1}.
 \eeq
 The fractional super group is the algebra $A^N_n(G)=A(G)\times\Lambda^N_n$
 where  $\Lambda^N_n$ is the algebra
 generated by  $\theta_\alpha$, $\alpha=1,\dots, N$, $\lambda$ with relations
 \beq\lb{dual1}
 \{\theta_{\alpha_1}, \theta_{\alpha_2},\dots,
 \theta_{\alpha_n}\}=0, \ \ \ \ \alpha_k\in 1,2,\dots, N
 \eeq
 and
 \beq\lb{dual2}
 \lambda\theta_\alpha =  q\theta_\alpha\lambda, \ \ \
 \lambda^n=1.
 \eeq
 Co-algebra operations and antipode in $A^N_n(G)$ depend
 on the structure constants $c^k_{ij}$, $a^j_{\alpha\beta}$ and
 $b^j_{\alpha_1\alpha_2,\dots,\alpha_n}$.
 As an example let us consider  the fractional super algebra
 \beq
 Q^n=P, \ \ \ [X, Q]=0.
 \eeq
 Since we have only one super element $Q$ the $S_n$ invariant form
 is equal up to  the multiple to $Q^n$.  The
 corresponding fractional group is  generated by $\theta$, $z$ and
 $\lambda$ such that
 \beq
 \theta^n=0, \ \ \ \lambda^n=1, \ \ \ \
 \lambda\theta=q\theta\lambda, \ \ \ q^n=1
 \eeq
 with $z$ being commutative with $\theta$ and $\lambda$. The  duality
 relations are
 \beq
 \langle Q, \theta \rangle =1, \ \ \ \langle X, z \rangle =1, \ \ \
 \langle K, \lambda \rangle =q.
 \eeq
 Using properties of duality relations we arrive at the following
 coalgebra structure
 \beq
 \Delta(\theta)=\theta\otimes 1 +\lambda\otimes\theta
 \eeq
 \beq
 \Delta(z)=z\otimes 1 +1\otimes z +\sum_{k=1}^{n-1}\fr{\lambda^{n-k}\theta^k\otimes\theta^{n-k}}
 {(q:q)_k(q:q)_{n-k}},
 \eeq
 \beq
 \varepsilon (\theta)=0, \ \ \ \varepsilon (\lambda )=1
 \eeq
 and
 \beq
 S(\theta)=-\lambda^{n-1}\theta, \ \  \ S(\lambda)=\lambda^{n-1}
 \eeq
 where
 \beq
 (q;q)_k=\prod_{j=1}^k (1-q^j)
 \eeq

 \vspace{5mm}
 \noindent
 {\bf 4. $S_3$ graded super algebras $sl(2)$}

 \vspace{2mm}
 \noindent
 From the commutation relations
 \beq\lb{s2a}
 [X_1, X_2]=X_3, \ \ \ [X_3, X_1]=2X_1, \ \ \ [X_3, X_2]=-2X_2
 \eeq
 for the algebra $sl(2)$ we read
 \beq
 c^3_{12}=1, \ \ \ c^1_{31}=2, \ \ \ c^2_{32}=-2.
 \eeq
 For given $N$ the matrix $a^j=\{a^j_{\alpha\beta}\}$ due to (\ref{rep1}) is an arbitrary
 $N$-dimensional representation of $sl(2)$.  The solution of (\ref{s3j1}) and (\ref{s3j2})
 for $b^j_{\alpha\beta\gamma}$ is fully determined by this
 representation. Since $b^j_{\alpha\beta\gamma}$ is symmetric in
 $\alpha,\ \beta$ and $\gamma$ through (\ref{s3al1}) the number of
 unknown coefficients for the $\mathrm{sl(2)}$ case is
 $N(N+1)(N+2)/2$. On the other hand, equation (\ref{s3j1}) which is
 symmetric in $\alpha,\ \beta,\ \gamma$ gives $3N(N+1)(N+2)/2$
 equations and equation (\ref{s3j2}) which is symmetric in
 $\alpha,\ \beta,\ \gamma,\ \sigma$ gives $N^{2}(N+1)(N+2)(N+3)/24$
 equations. Although the system seems overdetermined there are
 solutions some of which will be given below.
 We consider  $N=1, \ 2$ and $3$  fractional super  generalizations of $sl(2)$ at
 $n=3$, that is $q=e^{i\fr{2\pi}{3}}$.

 \vspace{2mm}
 \noindent
 {\bf A1. $N=1$ fractional super $sl(2)$}

 \vspace{1mm}
 \noindent
 We have one super generator $Q_1$ which can  transform  as scalar only. Therefore
 $a^j_{11}=0$.  Inserting it in the  relations (\ref{s3j1}) and  (\ref{s3j2}) we get
 $b^j_{111}=0$. These structure constants imply that the fractional super
 algebra $U_3^1(sl(2))$ is  the direct product  of the universal enveloping algebra $U(sl(2))$
 and the Hopf  algebra generated by $Q_1$  and $K$ satisfying the relations
 \beq
 KQ_1=qQ_1K, \ \ \  \ Q_1^3 =0, \ \ \ \ K^3=1
 \eeq
 and the co-algebra operations (\ref{com}) and (\ref{ant}).
 The fractional super group $A_1^3(SL(2))$ is the direct product of
 the Hopf algebras $A(SL(2))$ and  $\Lambda^1_3$. Recall that the Hopf algebra structure
 of polynomial algebra  $A(SL(2))$ is given by
 \beq\lb{z1}
 \Delta (x_{nm}) = \sum_{k=1}^2 x_{nk}\otimes x_{km}
 \eeq
 and
 \beq\lb{z11}
 S(x_{11})=x_{22}, \ \ S(x_{22})=x_{11}, \ \ S(x_{12})=-x_{12}, \ \
 S(x_{21})=-x_{21},
 \eeq
 where two by two matrix $x=\{x_{nm}\}$ representing $SL(2)$
 has  determinant $1$.
 The Hopf algebra structure of the
 algebra $\Lambda_3^1$ is given by the following formulas
 \beq
 \Delta(\theta_1)= \theta_1\otimes 1 + \lambda\otimes \theta_1, \ \ \
 \Delta(\lambda)=\lambda\otimes \lambda
 \eeq
 \beq
 S(\theta_1)=-\lambda^2\theta_1, \ \ \  S(\lambda)=\lambda^2.
 \eeq

 \vspace{2mm}
 \noindent
 {\bf A2. $N=2$ fractional super $sl(2)$}

 \vspace{1mm}
 \noindent
 For $N=2$ we have two possibilities. We can either  require
 generators $Q_1$, $Q_2$ to transform as scalars or as spinors.

 (i)  In the former case we have $a_{\alpha,\beta}^j=0$.
 From  the  relations (\ref{s3j1}) and  (\ref{s3j2}) we get
 $b^j_{\alpha\beta\gamma}=0$. The obtained structure constants
 imply that the fractional super algebra $U_3^2(sl(2))$ is the direct product
 of the universal enveloping algebra $U(sl(2))$ and the Hopf
 algebra generated by $Q_1$, $Q_2$ and $K$ satisfying the relations
 \beq
 KQ_\alpha=qQ_\alpha K, \ \ \  \ \{Q_\alpha,Q_\beta, Q_\gamma \}=0, \ \ \ \ K^3=1
 \eeq
 and the co-algebra operations (\ref{com}) and (\ref{ant}).
 The fractional super group $A_2^3(SL(2))$ is the direct product of
 the Hopf algebras $A(SL(2))$ and  $\Lambda^2_3$.
 The Hopf algebra structure of $\Lambda_3^2$ is given by the following formulas
 \beq
 \Delta(\theta_\alpha)= \theta_\alpha\otimes 1 + \lambda\otimes \theta_\alpha, \ \ \
 \Delta(\lambda)=\lambda\otimes \lambda
 \eeq
 \beq
 S(\theta_\alpha)=-\lambda^2\theta_\alpha, \ \ \  S(\lambda)=\lambda^2.
 \eeq

 (ii) Now let us assume that $Q_1$ and $Q_2$ transform as spinors under the action of $sl(2)$.
 We have
 \beq
 a^1 =
 \left(
 \begin{array}{cc}
 0 & 1 \\
 0 & 0
 \end{array}
 \right), \ \ \ \
 a^2 =
 \left(
 \begin{array}{cc}
 0 & 0 \\
 1 & 0
 \end{array}
 \right), \ \ \ \
 a^3 =
 \left(
 \begin{array}{cc}
 1 & 0 \\
 0 & -1
 \end{array}
 \right)
 \eeq
 Equation (\ref{s3j2}) gives 10 equations for 12 unknowns
 \beq\lb{sol1}
 b^1_{111}=b^1_{112}=b^2_{122}=b^2_{222}=b^3_{111}=b^3_{222}=0,
 \eeq
 \beq\lb{sol2}
 b^1_{122}=-\frac{1}{3}b^2_{111}=b^3_{112}, \ \
 b^1_{222}=-3b^2_{112}=3b^3_{122}
 \eeq
 Substituting these into (\ref{s3j1}) one finds that the only
 solution is  $b^j_{\alpha\beta\gamma}=0$.
 Thus we  obtained the following fractional super algebra
 \beq
 \{Q_\alpha,Q_\beta,Q_\gamma\}=0,
 \eeq
 \beq\lb{s2b}
 [Q_1, X_1]=Q_2, \ \ [Q_2, X_2]=Q_1, \ \ [Q_1, X_3]=Q_1, \ \  [Q_2,X_3]=-Q_2.
 \eeq
 Using the general construction given in the previous section one can define the fractional super
 group  $A^1_3(SL(2))$ corresponding to the above fractional super algebra.  $A^1_3(SL(2)$ is
 the algebra  generated by elements  $x_{nm}$, $\theta_n$, $n, m =1,2$  and  $\lambda$
 satisfying (\ref{dual1}), (\ref{dual2}) and $\det (x_{nm})=1$.
 The co-algebra operations and antipode can be shown to be given by (\ref{z1}),  (\ref{z11}) and
 \begin{eqnarray}
 \Delta( \theta_1) &=& \theta_2\otimes x_{21}
 +\theta_1\otimes x_{11}+\lambda\otimes\theta_1  \\
 \Delta( \theta_2) &=& \theta_2\otimes x_{22}
 +\theta_1\otimes x_{12}+\lambda\otimes\theta_2
 \end{eqnarray}
 and
 \beq
 S(\lambda)=\lambda^2, \ \
 S(\theta_1)=\lambda^2(x_{21}\theta_2-x_{22}\theta_1), \ \
 S(\theta_2)=\lambda^2(x_{12}\theta_1-x_{11}\theta_2).
 \eeq
 The duality relations are given by the formulas
 \begin{eqnarray}
 \langle X_3, x_{nn} \rangle = (-)^{n+1}, \ \
 \langle X_1, x_{12} \rangle = 1, \ \ \langle K, x_{nm} \rangle = \delta_{nm} \\
 \langle X_2, x_{21} \rangle = 1, \ \ \
 \langle Q_\alpha, \theta_\beta\rangle = \delta_{\alpha\beta}, \ \ \
 \langle K, \lambda \rangle = q.
 \end{eqnarray}

 \vspace{2mm}
 \noindent
 {\bf A3. $N=3$ fractional super $sl(2)$.}

 \vspace{1mm}
 \noindent
 We have three different superalgebras depending on the  choice of  $a^j$.

 (i) We take $a^j_{\alpha\beta}=0$. The relations (\ref{s3j1}) and
 (\ref{s3j2}) imply $b^j_{\alpha\beta\gamma=0}$. This case is
 similar with (i) of $A2$.

 (ii) We take the vector representation
 \beq
 a^1 =
 \left(
 \begin{array}{ccc}
 0 & 0 & 0\\
 \sqrt{2} & 0 & 0\\
 0 & \sqrt{2} & 0
 \end{array}
 \right), \
 a^2 =
 \left(
 \begin{array}{ccc}
 0 & \sqrt{2} & 0\\
 0 & 0 & \sqrt{2}\\
 0 & 0 & 0
 \end{array}
 \right), \
 a^3 =
 \left(
 \begin{array}{ccc}
 -2 & 0 & 0\\
 0 & 0 & 0\\
 0 & 0 & 2
 \end{array}
 \right).
 \eeq
 \noindent The substitution of

\beq\lb{cs1}
 a^1_{21}=a^1_{32}=a^2_{12}=a^2_{23}=\sqrt{2},\ \
a^3_{11}=-2,\ \ a^3_{33}=2
 \eeq

 \noindent into (\ref{s3j2}) gives

\begin{eqnarray}\lb{cs2}
b^1_{111} &=& 3\sqrt{2}b^3_{112}=-3b^2_{113}, \nonumber \\
b^1_{112} &=&\sqrt{2}b^3_{122}=-2b^2_{123}, \nonumber \\
b^1_{122}&=&\frac{\sqrt{2}}{3}b^3_{222}=-b^2_{223}, \\
b^1_{113}&=&2\sqrt{2}b^3_{123}=-b^2_{133}, \nonumber \\
b^1_{133}&=&\sqrt{2}b^3_{233}=-\frac{1}{3}b^2_{333}, \nonumber \\
b^1_{123}&=&\frac{\sqrt{2}}{2}b^3_{223}=-\frac{1}{2}b^2_{233}\nonumber
\end{eqnarray}

\noindent and the remaining 12 parameters
$b^j_{\alpha\beta\gamma}$ are zero. The substitution of
(\ref{cs2}) into (\ref{s3j1}) gives

\beq\lb{cs3}
b^1_{113}=-2b^1_{122}=-b^2_{133}=2b^2_{223}=2\sqrt{2}b^3_{123}=-\frac{2\sqrt{2}}{3}b^3_{222}
\eeq

\noindent and all other $b^j_{\alpha\beta\gamma}$ are zero. Hence
we have a unique extension for the vector representation of
$\mathrm{sl(2)}$. Equations satisfied by
$b^j_{\alpha\beta\gamma}$, namely (\ref{s3j1}) and (\ref{s3j2}),
are invariant under rescaling $b^j_{\alpha\beta\gamma}\rightarrow
kb^j_{\alpha\beta\gamma}$ where $k$ is any nonzero constant. The
choice of this nonzero constant results only in a rescaling of the
generators $Q_{\alpha}$. Just for the sake of simplicity we choose
$Q_2^3=X_3$, i.e., $b^3_{222}=6$. Then the fractional
supersymmetric extension of $\mathrm{sl(2)}$ reads

\beq\lb{cst3}
 [Q_1, X_2]=\sqrt{2}Q_2, \  [Q_1, X_3]=-2Q_1, \  [Q_2, X_1]=\sqrt{2}Q_1,
\eeq

\beq\lb{cst4} [Q_2, X_2]=\sqrt{2}Q_3, \  [Q_3, X_1]=\sqrt{2}Q_2, \
[Q_3,X_3]=2Q_3 \eeq

\noindent and

\beq \lb{cst5}
 \{Q_1,Q_1,Q_3\}=-4\sqrt{2}X_1, \
\{Q_1,Q_2,Q_2\}=2\sqrt{2}X_1,   \ \{Q_1,Q_2,Q_3\}=-2X_3, \eeq

\beq\lb{cst6}
 \{Q_1,Q_3,Q_3\}=-4\sqrt{2}X_2, \
\{Q_2,Q_2,Q_2\}=6X_3, \ \{Q_2,Q_2,Q_3\}=-2\sqrt{2}X_2. \eeq

\noindent Notice also that all $b^j_{\alpha\beta\gamma}=0$ is
always a solution of (\ref{s3j1}) and (\ref{s3j2}).

 (iii) Assume that two of fractional  super generators  $Q_1$, $Q_2$ and $Q_3$ transform
 as spinors and the remaining one transforms as scalar, that is
 \beq
 a^1 =
 \left(
 \begin{array}{ccc}
 0 & 1 & 0\\
 0 & 0 & 0\\
 0 & 0 & 0
 \end{array}
 \right), \
 a^2 =
 \left(
 \begin{array}{ccc}
 0 & 0 & 0\\
 1 & 0 & 0\\
 0 & 0 & 0
 \end{array}
 \right), \
 a^3 =
 \left(
 \begin{array}{ccc}
 1 & 0 & 0\\
 0 & -1 & 0\\
 0 & 0 & 0
 \end{array}
 \right).
 \eeq
 The conditions (\ref{s3j1}) and (\ref{s3j2}) imply
 \beq
 b^1_{223}=-b^2_{111}=2b^3_{123}
 \eeq
 with all other structure coefficients $b^j_{\alpha\beta\gamma}$
 being zero. Choosing $b^1_{223}=1$ we
 get  the fractional super  algebra  given
 by (\ref{s2b}) and
 \beq\lb{s2c}
 \{Q_1,Q_1,Q_3\}=-X_2, \ \ \{Q_2,Q_2,Q_3\}=X_1, \ \ \
 \{Q_1,Q_2,Q_3\}=\fr{1}{2}X_3
 \eeq
 with all other relations being zero.

 Before closing this section we discuss realization of fractional super algebras  by
 "differential operators" in some linear spaces. Recall that for realization of super algebras
 one uses super derivatives which act on superspaces. Let $F(M)$
 be an algebra of functions on a manifold $M$. For fixed grading $n$ and the number $N$
 of  "grassmannian" variables  a fractional superspace is defined to be the direct product
 algebra  $F(M)\times \Lambda_n^N$. We define  fractional
 derivatives $D_{\theta_\alpha}$ by the formulas
 \beq
 D_{\theta_\alpha}  \theta_\beta  = \delta_{\alpha\beta}, \ \ \ \
 D_{\theta_\alpha} (ab)=D_{\theta_\alpha}(a) b + k(a)
 D_{\theta_\alpha}(b),
 \eeq
 where $a, \ b\in\Lambda_n^N$ and
 \beq
 k (\theta_\alpha)=q\theta_\alpha, \ \ \ \ k(ab)=k(a)k(b).
 \eeq
 Note that  $D_{\theta_\alpha}(f)=0$ and $k(f)=f$ if $f\in F(M)$.
 One can verify that these derivatives satisfy the relations
 \beq
 \sum_{\alpha_1\dots\alpha_n\in S_n}D_{\theta_{\alpha_1}}\cdots
 D_{\theta_{\alpha_n}}=0
 \eeq
 By means of fractional derivatives and  superspaces
 defined above one  can construct a realization of a fractional
 superalgebras.  For example the formulas
 \beq\lb{E}
 X_1 = -z^2\fr{d}{dz}-zL, \ \ \ X_2 = \fr{d}{dz}, \ \ \
 X_3 = 2z\fr{d}{dz}+ L
 \eeq
 \beq
 Q_1 = D_\theta, \ \  \ Q_2 = -zD_\theta, \ \ \  Q_3 =\fr{q}{2}\theta^2
 \fr{d}{dz}, \ \ \ K=q^L,
 \eeq
 where  $q=e^{i\fr{2\pi}{3}}$ and
 \beq
 L=-q(2\theta^2D_\theta^2+D_\theta\theta^2D_\theta)
 \eeq
 define representation of the fractional algebra (iii) in the linear
 space $A(C)\times \Lambda_3^1$, where $A(C)$ is the algebra of
 polynomials of the complex variable $z$. Indeed using
 \beq\lb{gr}
 L \theta^k = k\theta^k
 \eeq
 and
 \beq
 \theta^2 D^2_\theta +D^2_\theta\theta^2 +  D_\theta\theta^2D_\theta =-q^2
 \eeq
 one can easily  verify  the relations  (\ref{s2a}), (\ref{s2b}) and
 (\ref{s2c}).


\newpage

 \begin{thebibliography}{99}
 \bibitem{rau1} M. Raush deTraunbenberg and M.J. Slupinski,
              {\it J. Math. Phys.}, {\bf 41}, 4556 (2000).
 \bibitem{rau2} M. Raush de Traunbenberg,
 Fractional sypersymmetry and Lie algebras, arXiv:hep-th/0007150.
 \bibitem{dewitt} B. De Witt,
  "Supermanifols", Cambridge Univ. Press, Cambridge (1984).
 \bibitem{kos} B. Kostant, "Graded manifolds, graded Lie theory and
 prequantization," in Lecture Notes in Mathematics, Spriger
 Verlag, New York, vol. 570, 177 (1977).
 \bibitem{bat} M. Batchelor,
  {\it Tran. Am. Math. Soc.}, {\bf 258}, 257 (1980).
  \bibitem{vil1} N. Ya. Vilenkin and A. U. Klimyk,
 Representations of Lie Groups and Special Functions, vol 3, Kluwer
 Academic Press, The Netherland (1991).
 \bibitem{ahn} C. Ahn, D. Bernard and A. Leclair,
              {\it Nucl. Phys. }, {\bf B346}, 409 (1990).
 \bibitem{rau} M. Raush de Traunbenberg and M.J. Slupinski,
              {\it Mod. Phys. Lett. A}, {\bf 39}, 3051 (1997).
 \bibitem{ker} R. Kerner,
              {\it J. Math. Phys.}, {\bf 33}, 403 (1992).
 \bibitem{vik} V. Abramov, R. Kerner and B. Le Roy
              {\it J. Math. Phys.}, {\bf 38}, 1650 (1997).
 \bibitem{azc} J.A. de Azcarraga and M.J. Macfarlane,
              {\it J. Math. Phys.}, {\bf 37}, 1115 (1996).
 \bibitem{dur} S. Durand,
              {\it  Mod. Phys. Lett.}, {\bf A7}, 2905 (1992).
 \bibitem{ahm1} H. Ahmedov and \"{O}.F. Dayi,
              {\it J. Phys. A}, {\bf 32}, 1895 (1999);
              {\it J. Phys. A}, {\bf 32}, 6247 (1999);
  {\it Mod. Phys. Lett. A}, {\bf 15}, no. 29, 1801 (2000);
 H. Ahmedov,  {\it Turkish J. Phys.}, {\bf 24}, no. 3, 175 (2000)
 \bibitem{vil} N. Ya. Vilenkin and A. U. Klimyk,
 Representations of Lie Groups and Special Functions, vol 3, Kluwer
 Academic Press, The Netherland (1991).

 \end{thebibliography}
 \end{document}